\pdfoutput=1

\documentclass[letterpaper, conference, 10pt]{IEEEtran}

\usepackage[cmex10]{amsmath}
\usepackage{amsthm}

\usepackage{amssymb}
\usepackage{color}
\usepackage{graphicx}
\usepackage{ifthen}
\usepackage{mathrsfs} 
\usepackage{mathscinet} 
\usepackage{mathtools}
\usepackage{rotating}
\usepackage{stmaryrd}
\usepackage{cite}

\mathtoolsset{%
   mathic=true
}%

\newcommand*{\Z}{\mathbb{Z}}
\newcommand*{\Q}{\mathbb{Q}}
\newcommand*{\R}{\mathbb{R}}

\newcommand*{\verts}[1]{\left\lvert #1 \right\rvert}
\newcommand*{\braces}[1]{\left\lbrace #1 \right\rbrace}

\newcommand*{\floors}[1]{\left\lfloor #1 \right\rfloor}

\newcommand*{\card}{\verts}
\newcommand*{\setof}{\braces}

\newcommand*{\deftobe}{\mathrel{\coloneqq}}

\newcommand*{\maps}{\colon}
\newcommand*{\st}{\,:\,}

\renewcommand*{\epsilon}{\varepsilon}
\renewcommand*{\phi}{\varphi}

\newtheorem{thm}{Theorem}[section]

\newtheorem{lem}[thm]{Lemma}

\newtheorem{prop}[thm]{Proposition}

\theoremstyle{definition}

\newtheorem{defn}[thm]{Definition}

\theoremstyle{remark}

\newcommand*{\pZ}{\mathrm{p}\Z}
\newcommand*{\D}{\mathcal{D}}
\newcommand*{\SL}{\mathsf{SL}}
\newcommand*{\GL}{\mathsf{GL}}

\DeclareMathOperator*{\latlen}{len}  
\newcommand*{\bnd}{\mathfrak{b}}  
\newcommand*{\intr}{\mathfrak{I}} 
\newcommand*{\A}{\mathfrak{A}}    
\newcommand*{\inhl}[1]{\widetilde{#1}}   
\newcommand*{\tran}{\mathrm{t}}  
\renewcommand*{\card}[1]{\verts{#1}}
\DeclareMathOperator{\conv}{Conv}

\newcommand*{\Ehr}{\mathfrak{L}}
\DeclareMathOperator{\lcm}{lcm}

\title{%
   {\huge \textbf{Coefficient functions of the Ehrhart
   quasi-polynomials of rational polygons}}
}%

\author{%
   \IEEEauthorblockN{{\large \textbf{%
      Tyrrell B. McAllister%
   }}}
   \IEEEauthorblockA{{\large%
      Department of Mathematics and Computer Science,
      Eindhoven University of Technology
   } \\ {\large
      Eindhoven,
      The Netherlands
   }}
   \thanks{%
      Research supported by the Netherlands Organisation for
      Scientific Research (NWO) Mathematics Cluster DIAMANT.%
   }
}%

\date{Draft}

\begin{document}
   \maketitle

   \begin{abstract}
      In 1976, P.~R.~Scott characterized the Ehrhart polynomials
      of convex integral polygons.  We study the same question for
      Ehrhart polynomials and quasi-polynomials of
      \emph{non}-integral convex polygons.
      Define a \emph{pseudo-integral polygon}, or \emph{PIP}, to
      be a convex rational polygon whose Ehrhart quasi-polynomial
      is a polynomial.  The numbers of lattice points on the
      interior and on the boundary of a PIP determine its Ehrhart
      polynomial.
      We show that, unlike the integral case, there exist PIPs
      with $b=1$ or $b=2$ boundary points and an arbitrary number
      $I \ge 1$ of interior points.  However, the question of
      whether a PIP must satisfy Scott's inequality $b \le 2I + 7$
      when $I \ge 1$ remains open.
      Turning to the case in which the Ehrhart quasi-polynomial
      has nontrivial quasi-period, we determine the possible
      minimal periods that the coefficient functions of the
      Ehrhart quasi-polynomial of a rational polygon may have.
   \end{abstract}
   

   \section{%
      Introduction%
   }%
   
   We take a \emph{rational polygon} $P \subset \R^{2}$ to be the
   convex hull of finitely many rational points, not all contained
   in a line.  In particular, all of our polygons are convex.
   Given a positive integer $n$, let $nP \deftobe \{nx \in \R^{2}
   \st x \in P\}$ be the dilation of $P$ by $n$.  The
   $2$-dimensional case of a well-known result due to E.~Ehrhart
   \cite{Ehr62} states that the number $\card{nP \cap \Z^{2}}$ of
   integer lattice points in $nP$ is a degree-$2$ quasi-polynomial
   function of $n$ with rational coefficients.  That is, there
   exist periodic functions $c_{P,i} \maps \Z \to \Q$, $i = 0, 1,
   2$, such that, for all positive integers $n$,
   \begin{align*}
      \Ehr_{P}(n)
      &\deftobe
	 c_{P,2}(n) n^{2} 
	 +  c_{P,1}(n) n
	 +  c_{P,0}(n) \\
	 &= \card{nP \cap \Z^{2}}.
   \end{align*}
   We call $\Ehr_{P}$ the \emph{Ehrhart quasi-polynomial} of $P$.
   We say that $P$ has \emph{period sequence} $(s_{2}, s_{1},
   s_{0})$ if the minimum period of the coefficient function
   $c_{P,i}$ is $s_{i}$ for $i = 0, 1, 2$.  The
   \emph{quasi-period} of $\Ehr_{P}$ (or of $P$) is $\lcm
   \setof{s_{0}, s_{1}, s_{d}}$.  We refer the reader to
   \cite{BR07} for a thorough introduction to the theory of
   Ehrhart quasi-polynomials.
   
   Our goal in this note is to examine the properties and possible
   values of the coefficient functions $c_{P,i}$.  The leading
   coefficient $c_{P,2}$ is always the area $\A_{P}$ of $P$.
   Furthermore, when $P$ is an integral polygon (meaning that its
   vertices are all integer lattice points), $\Ehr_{P}$ is simply
   a polynomial with $c_{P,0}=1$ and $c_{P,1} = \frac{1}{2}
   \bnd_{P}$, where $\bnd_{P}$ is the number of integer lattice
   points on the boundary of $P$.  Now, Pick's formula determines
   $\A_{P}$ in terms of $\bnd_{P}$ and the number $\intr_{P}$ of
   integer lattice points in the interior of $P$.  Hence,
   characterizing the Ehrhart polynomials of integral polygons
   amounts to determining the possible numbers of integer lattice
   points in their interiors and on their boundaries.  This was
   accomplished by P.~R.~Scott \cite{Sco76} in 1976:
   
   \begin{thm}[P.~R.~Scott \cite{Sco76}]
      Given non-negative integers $I$ and $b$, $(I, b) =
      (\intr_{P}, \bnd_{P})$ for some integral polygon $P$ if and
      only if $b \ge 3$ and either $I = 0$, $(I, b) = (1, 9) $, or
      $b \le 2 I + 6$.
   \end{thm}
   
   However, not all Ehrhart polynomials of polygons come from
   \emph{integral} polygons.  Hence, the complete characterization
   of Ehrhart polynomials of rational polygons, including the
   non-integral ones, remains open.  To this end, we define a
   \emph{pseudo-integral polygon}, or \emph{PIP}, to be a rational
   polygon with quasi-period~1.  That is, PIPs are those polygons
   that share with integral polygons the property of having a
   polynomial Ehrhart quasi-polynomial.  Like integral polygons,
   PIPs must satisfy Pick's Theorem \cite[Theorem 3.1]{MW05}, so,
   again, the problem reduces to finding the possible values of
   $\intr_{P}$ and $\bnd_{P}$.  In Section
   \ref{sec:NonintegralPIPs}, we construct PIPs with $\bnd_{P} \in
   \setof{1,2}$ and $\intr_{P}$ an arbitrary positive integer.
   This construction therefore yields an infinite family of
   Ehrhart polynomials that are not the Ehrhart polynomial of any
   integral polygon.
   
   In Section \ref{sec:PeriodsOfCoefficients}, we consider the
   case where $P$ is not a PIP\@.  Determining all possible
   coefficient functions $c_{P,i}$ seems out of reach at this
   time.  However, one interesting question that we will answer
   here is, what are the possible period sequences $(s_{2}, s_{1},
   s_{0})$?  P.~McMullen showed that $s_{i}$ is bounded by the
   so-called \emph{$i$-index} of $P$ \cite{McM78}.  We state his
   result here in the full generality of $d$-dimensional
   polytopes:
   \begin{thm}[McMullen \protect{\cite[Theorem 6]{McM78}}]
   \label{thm:McMullenBound}
      Given a $d$-dimensional polytope $P$ and $i \in \setof {0,
      \dotsc, d}$, define the \mbox{\emph{$i$-index}} of $P$ to be
      the least positive integer $p_{i}$ such that all the
      $i$-dimensional faces of $p_{i}P$ contain integer lattice
      points in their affine span.  Then the period $s_{i}$ of the
      $i$th coefficient of $\Ehr_{P}$ divides $p_{i}$.  In
      particular, $s_{i} \le p_{i}$.
   \end{thm}
   Observe that, by definition, $p_{d} \mid p_{d-1} \mid \dotsb
   \mid p_{0}$.  Conversely, Beck, Sam, and Woods \cite{BSW08}
   have shown that, given any positive integers $p_{d} \mid
   p_{d-1} \mid \dotsb \mid p_{0}$, there exists a polytope with
   $i$-index $p_{i}$ for $0 \le i \le d$.  Moreover, McMullen's
   bounds on the $s_{i}$'s are tight for this polytope: $s_{i} =
   p_{i}$.
   
   Thus we have that $s_{i}$ is bounded by the $i$-index, and this
   bound is tight in some cases.  Furthermore, the $i$-index
   weakly increases as $i$ decreases.  Seeing this, one might hope
   that the $s_{i}$'s themselves are also required to satisfy some
   constraints.  However, in Section
   \ref{sec:PeriodsOfCoefficients}, we show that, in the case of
   polygons, $s_{0}$ and $s_{1}$ may take on arbitrary values.

   \section{%
      Piecewise skew unimodular transformations
   }%

   Since we will be exploring the possible Ehrhart
   quasi-polynomials of polygons, it will be useful to have a
   geometric means of constructing polygons while controlling
   their Ehrhart quasi-polynomials.  The main tool that we will
   use are piecewise affine unimodular transformations.  Following
   \cite{Gre93}, we call these \emph{$\pZ$-homeomorphisms}.
   \begin{defn}
      Given $U, V \subset \R^{2}$ and a finite set
      $\setof{\ell_{i}}$ of lines in the plane, let
      $\setof{C_{j}}$ be the set of connected components of $U
      \setminus \bigcup_{i} \ell_{i}$.  Then a homeomorphism
      $f\maps U \to V$ is a \emph{$\pZ$-homeomorphism} if, for
      each component $C_{j}$, $f \vert_{C_{j}}$ is the restriction
      to $C_{j}$ of an element of $\GL_{2}(\Z) \ltimes \Z^{2}$.
   \end{defn}
   The key property of $\pZ$-homeomorphisms is that they preserve
   the lattice and so preserve Ehrhart quasi-polynomials.
   
   In particular, we will be using $\pZ$-homeomorphisms that act
   as skew transformations on each component of their domains.
   Given a rational vector $r \in \Q^{2}$, let $r_{p}$ be the
   generator of the semigroup $(\R_{\ge 0}r) \cap \Z^{2}$, and
   define the \emph{lattice length} $\latlen(r) \in \Q$ of $r$ by
   $\latlen(r) r_{p} = r$.  Thus, if $r = (\tfrac{a}{b},
   \tfrac{c}{d})$, where the fractions are reduced, then we have
   that $\latlen(r) = \gcd(a, c) / \lcm(b, d)$.  Define the skew
   unimodular transformation $U_{r} \in \SL_{2}(\Z)$ by
   \begin{equation*}
      U_{r} (x) 
      =  x 
         +  \frac{1}{\latlen(r)^{2}}
            \det(r, x) r,
   \end{equation*}
   where $\det(r,x)$ is the determinant of the matrix whose
   columns are $r$ and $x$ (in that order).  Equivalently, let $S$
   be the subgroup of skew transformations in $\SL_{2}(\Z)$ that
   fix $r$.  Then $U_{r}$ is the generator of $S$ that translates
   $v$ parallel (resp.\ anti-parallel) to $r$ if the angle between
   $r$ and $v$ is less than (resp.\ greater than) 180 degrees
   measured counterclockwise.
   
   Define the piecewise unimodular transformations $U_{r}^{+}$ and
   $U_{r}^{-}$ by
   \begin{equation*}
      U_{r}^{+}(x)
      =  \begin{dcases*}
            U_{r}(x) & if $\det(r,x) \geq 0$,  \\
            x & else,
         \end{dcases*}
   \end{equation*}
   and
   \begin{equation*}
      U_{r}^{-}
      =  (U_{-r}^{+})^{-1}.
   \end{equation*}
   Finally, given a lattice point $u \in \Z^{2}$ and a rational
   point $w \in \Q^{2}$, let $U_{uw}^{+}$ and $U_{uw}^{-}$ be the
   \emph{affine} piecewise unimodular transformations defined by
   \begin{align*}
      U_{uw}^{+}(v)
      &= U_{w - u}^{+}(v - u) + u,  \\
      U_{uw}^{-}(v) 
      &= U_{w - u}^{-}(v - u) + u.
   \end{align*}

   \section{%
       Constructing nonintegral PIPs%
   }%
   \label{sec:NonintegralPIPs}

   \begin{thm}
      There does not exist a 2-dimensional PIP $P$ with $\bnd_{P}
      = 0$ or with $(\intr_{P}, \bnd_{P}) \in \{ (0,1), (0, 2)
      \}$.  However, for all integers $I \ge 1$ and $b \in \{1,
      2\}$, there exists a PIP $P$ with $(\intr_{P}, \bnd_{P}) =
      (I, b)$.
   \end{thm}
   
\begin{figure}[tbp]
      \begin{center}
      \includegraphics{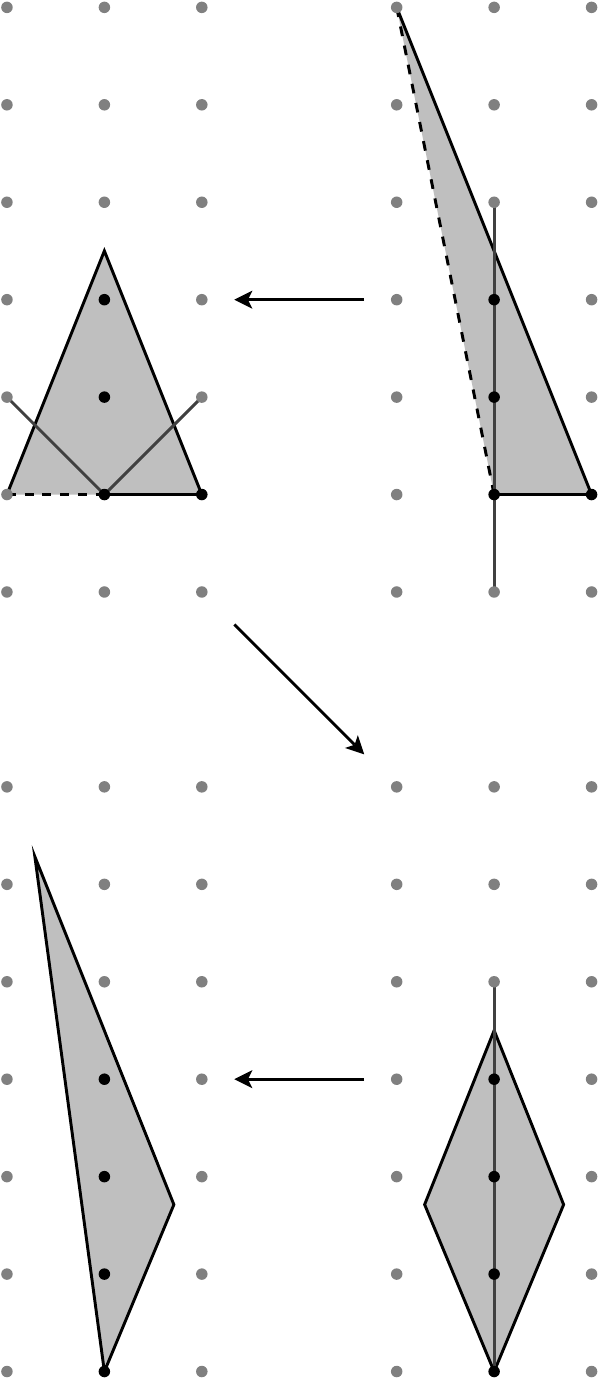}
      \caption{%
	 The construction of a PIP with one boundary point and an
	 arbitrary number $I$ of interior points in the case $I =
	 3$.
      }%
      \label{fig:OneVertex}
      \end{center}
   \end{figure}

   \begin{IEEEproof}
      In \cite[Theorem 3.1]{MW05}, it was shown that if $P$ is a
      PIP, then $\bnd_{nP} = n \bnd_{P}$ for $n \in \Z_{>0}$.  If
      $\bnd_{P} = 0$, this implies that $\bnd_{nP} = 0$ for all $n
      \in \Z_{>0}$, which is impossible because, for example, some
      multiple $\D P$ of $P$ is integral.  Hence, $\bnd_{P} \ge 1$
      when $P$ is a PIP.
       
      It was also shown in \cite{MW05} that PIPs satisfy Pick's
      theorem.  Hence, we must have that $\A_{P} = \intr_{P} +
      \frac{1}{2} \bnd_{P} - 1$.  But if $\intr_{P} = 0$ and
      $\bnd_{P} \in \{1, 2\}$, this yields an area less than or
      equal to 0.  Since we are not considering polygons contained
      in a line, this is impossible.
      
      Therefore, if $\bnd_{P} < 3$, we must have $\bnd_{P} \in
      \{1, 2\}$ and $\intr_{P} \ge 1$.  Now let integers $b \in
      \setof{1,2}$ and $I \ge 1$ be given.  We construct a PIP $P$
      with $(\intr_{P}, \bnd_{P}) = (I,b)$.
      
      If $b = 2$, consider the triangle \[T = \conv
      \setof{(0,0)^{\tran}, (I+1, 0)^{\tran}, (1, 1 -
      \tfrac{1}{I+1})^{\tran}}.\] It was proved in \cite{MW05}
      that $T$ is a PIP\@.  Let $P$ be the union of $T$ and its
      reflection about the $x$-axis.  Then $\intr_{P} = I$ and
      $\bnd_{P} = 2$.  Moreover, $\Ehr_{P}(n) = 2\Ehr_{T}(n) - I -
      2$ (correcting for points double-counted on the $x$-axis),
      so $P$ is also a PIP.
      
      If $b = 1$, consider the ``semi-open'' triangle 
      \begin{align*}
          T_{1}
	  &= 
	     \conv\setof{(0,0)^{\tran}, (1, 2I - 1)^{\tran}, (-1, 0)} \\
          &\quad \setminus
	     \bigl( (0,0)^{\tran}, (1, 2I - 
	     1)^{\tran} \bigr].
      \end{align*}
      (The upper left of Figure \ref{fig:OneVertex} depicts the
      case with $I = 3$.)  The Ehrhart quasi-polynomial of $T$ is
      evidently a signed sum of Ehrhart polynomials, so it also is
      a polynomial.  We will apply a succession of
      $\pZ$-homeomorphisms to $T$ to produce a convex rational
      polygon without changing the Ehrhart polynomial.  (The gray
      line-segments in Figure \ref{fig:OneVertex} indicate the
      lines that will be fixed by our skew transformations.)
      
      Let $T_{2} = (U_{(0,-1)^{\tran}}^{+})^{2I - 1}(T_{1})$.
      Hence, $T_{2} = \conv \setof{(1,0)^{\tran}, (0, I -
      1/2)^{\tran}, (-1,0)^{\tran}} \setminus \bigl((0,0)^{\tran},
      (1,0)^{\tran}\bigr]$.  (See Figure \ref{fig:OneVertex},
      upper right.)
      
      Now act upon the triangle below the line spanned by
      $(-1,-1)$ (resp.\ $(1,-1)$), with $U_{(-1,-1)^{\tran}}^{+}$
      (resp\ $U_{(1,-1)^{\tran}}^{-}$).  The result is
      \begin{multline*}
         T_{3}
         = 
         \conv \left\lbrace
            \begin{pmatrix}
               0  \\
               -1
            \end{pmatrix},\;
            \begin{pmatrix}
               \dfrac{2I - 1}{2I + 1}  \\
               \\
               2I \dfrac{2I - 1}{2I + 1}
            \end{pmatrix}, \right. \\
            \left.
            \begin{pmatrix}
               -\dfrac{2I - 1}{2I + 1}  \\
               \\
               2I \dfrac{2I - 1}{2I + 1}
            \end{pmatrix},\;
            \begin{pmatrix}
               0  \\
               I - 1/2
            \end{pmatrix}
         \right \rbrace,
      \end{multline*}
      (see Figure \ref{fig:OneVertex}, lower left).  At this
      point, we have a convex rational polygon with the desired
      number of interior and boundary lattice points, so the claim
      is proved.  However, it might be noted that we can get a
      triangle by letting $P = (U_{(0,1)^{\tran}}^{-})^{^{2I -
      1}}(T_{3})$, yielding
      \begin{equation*}
         P
         =  \conv
            \setof{
               \begin{pmatrix}
                  0  \\
                  -1
               \end{pmatrix},\;
	       \begin{pmatrix}
		  \dfrac{2I - 1}{2I + 1}  \\
		  \\
		 \dfrac{2I - 1}{2I + 1}
	       \end{pmatrix},\;
	       \begin{pmatrix}
		  -\dfrac{2I - 1}{2I + 1}  \\
		  \\
		  2I \dfrac{2I - 1}{2I + 1}
	       \end{pmatrix}
             }.
       \end{equation*}
   \end{IEEEproof}
   
   A proof of, or counterexample to, Scott's inequality for
   nonintegral PIPs eludes us.  However, it is easy to show that
   any counterexample $P$ cannot contain a lattice point in the
   interior of its integral hull $\inhl{P} \deftobe \conv (P \cap
   \Z^{2})$.
   
   \begin{prop}
      If $P$ is a polygon whose integral hull contains a lattice
      point in its interior, then $P$ obeys Scott's
      inequality---that is, $\bnd_{P} \le 2 \intr_{P} + 6$ unless
      $(\intr_{P}, \bnd_{P}) = (1, 9)$.
   \end{prop}
   
   \begin{IEEEproof}
      We are given that $\intr_{\inhl{P}} \ge 1$.  Note that
      $\bnd_{\inhl{P}} \ge \bnd_{P}$ and $\intr_{\inhl{P}} \le
      \intr_{P}$.  Since $\inhl{P}$ is an integral polygon, it
      obeys Scott's inequalities: either $\bnd_{\inhl{P}} \le 2
      \intr_{\inhl{P}} + 6$ or $(\intr_{\inhl{P}},
      \bnd_{\inhl{P}}) = (1, 9)$.  In the former case, we have
      $%
         \bnd_{P}
         \le \bnd_{\inhl{P}}
         \le 2 \intr_{\inhl{P}} + 6
         \le 2 \intr_{P} + 6
      $.
      In the latter case, we similarly have $\bnd_{P} \le 9$ and
      $1 \le \intr_{P}$, so either $\intr_{P} = 1$ or $\bnd_{P}
      \le 2 \intr_{P} + 6$.
   \end{IEEEproof}

   \section{%
      Periods of coefficients of Ehrhart Quasi-polynomials%
   }%
   \label{sec:PeriodsOfCoefficients}

   If $P$ is a rational polygon, then the coefficient of the
   leading term of $\Ehr_{P}$ is the area of $P$, so the first
   term in the period sequence of $P$ is $1$.  However, we show
   below that no constraints apply to the remaining terms in the
   period sequence:
   \begin{thm}
   \label{thm:CoefficientPeriodSequences}
      Given positive integers $s$ and $t$, there exists a polygon
      $P$ with period sequence $(1, s, t)$.
   \end{thm}
  
   Before proceeding to the proof, we will need some elementary
   properties of the coefficients of certain Ehrhart
   quasi-polynomials.
   
   Fix a positive integer $s$, and let $\ell$ be the line segment
   $[0, \tfrac{1}{s}]$.  Then we have that $\Ehr_{\ell}(n) =
   \tfrac{1}{s} n + c_{\ell,0}(n)$, where the ``constant''
   coefficient function $c_{\ell,0}(n) = \floors{n/s} - n/s + 1$
   has minimum period $s$.  Note also that the half-open interval
   $h = (\tfrac{1}{s}, 1]$ satisfies $\Ehr_{\ell} + \Ehr_{h} =
   \Ehr_{[0,1]}$.  In particular, we have that
   \begin{equation}
      c_{\ell,0} + c_{h,0} = 1.
      \label{eq:ConstantCoeffofLineSegment}
   \end{equation}

   Given a positive integer $m$, it is straightforward to compute
   that the Ehrhart quasi-polynomial of the rectangle $\ell \times
   [0, m]$ is given by
   \begin{equation*}
      \Ehr_{\ell \times [0, m]}(n)
      =  \tfrac{m}{s}n^{2}
         +  \big(
               m c_{\ell,0} (n)
               +  \tfrac{1}{s}
            \big)
            n
         +  c_{\ell,0} (n).
   \end{equation*}
   In particular, the ``linear'' coefficient function has minimum
   period $s $, and the ``constant'' coefficient function is
   identical to that of $\Ehr_{\ell}$.  More strongly, we have the
   following:
   \begin{lem}
   \label{lem:CoeffsofRectangleUnionIntegral}
      Suppose that a polygon $P$ is the union of $\ell \times [0,
      m]$ and an integral polygon $P'$ such that $P' \cap (\ell
      \times [0, m])$ is a lattice segment.  Then $c_{P,1}$ has
      minimum period $s$ and $c_{P,0} = c_{\ell,0}$.
   \end{lem}

   With these elementary facts in hand, we can now prove Theorem
   \ref{thm:CoefficientPeriodSequences}.
   
   \begin{IEEEproof}[%
      Proof of Theorem \ref{thm:CoefficientPeriodSequences}%
   ]%
      Any integral polygon has period sequence $(1,1,1)$, so we
      may suppose that either $s \ge 2$ or $t \ge 2$.  Our
      strategy is to construct a polygon $H$ with period sequence
      $(1, s, 1)$ and a triangle $Q$ with period sequence
      $(1,1,t)$.  We will then be able to construct a polygon with
      period sequence $(1,s,t)$ for $s, t \ge 2$ by gluing $H$ and
      $Q$ along an integral edge.
      
      We begin by constructing a polygon with period sequence $(1,
      s, 1)$ for an arbitrary integer $s \ge 2$.  Define $H$ to be
      the heptagon with vertices
      \begin{align*}
	 t_{1} 
	 &= \big(
	       -\tfrac{1}{s},\, s(s-1) + 1
	    \big)^{\tran}, \\
	 t_{2}
	 &= \big(
	       -\tfrac{1}{s},\, -s(s-1) - 1
	    \big)^{\tran}, \\
	 u_{1}
	 &= \big(
	       0,\, s(s-1) + 1
	    \big)^{\tran}, \\
	 u_{2}
	 &= \big(
	       0,\, -s(s-1) - 1
	    \big)^{\tran}, \\
	 v_{1}
	 &= \big(
	       1,\, s(s-1)
	    \big)^{\tran}, \\
	 v_{2}
	 &= \big(
	       1,\, -s(s-1)
	    \big)^{\tran}, \\
	 w 
	 &= \big(
	       s - 1 + \tfrac{1}{s},\, 0
	    \big)^{\tran}.
      \end{align*}
      \begin{figure*}
         \centering
         \includegraphics{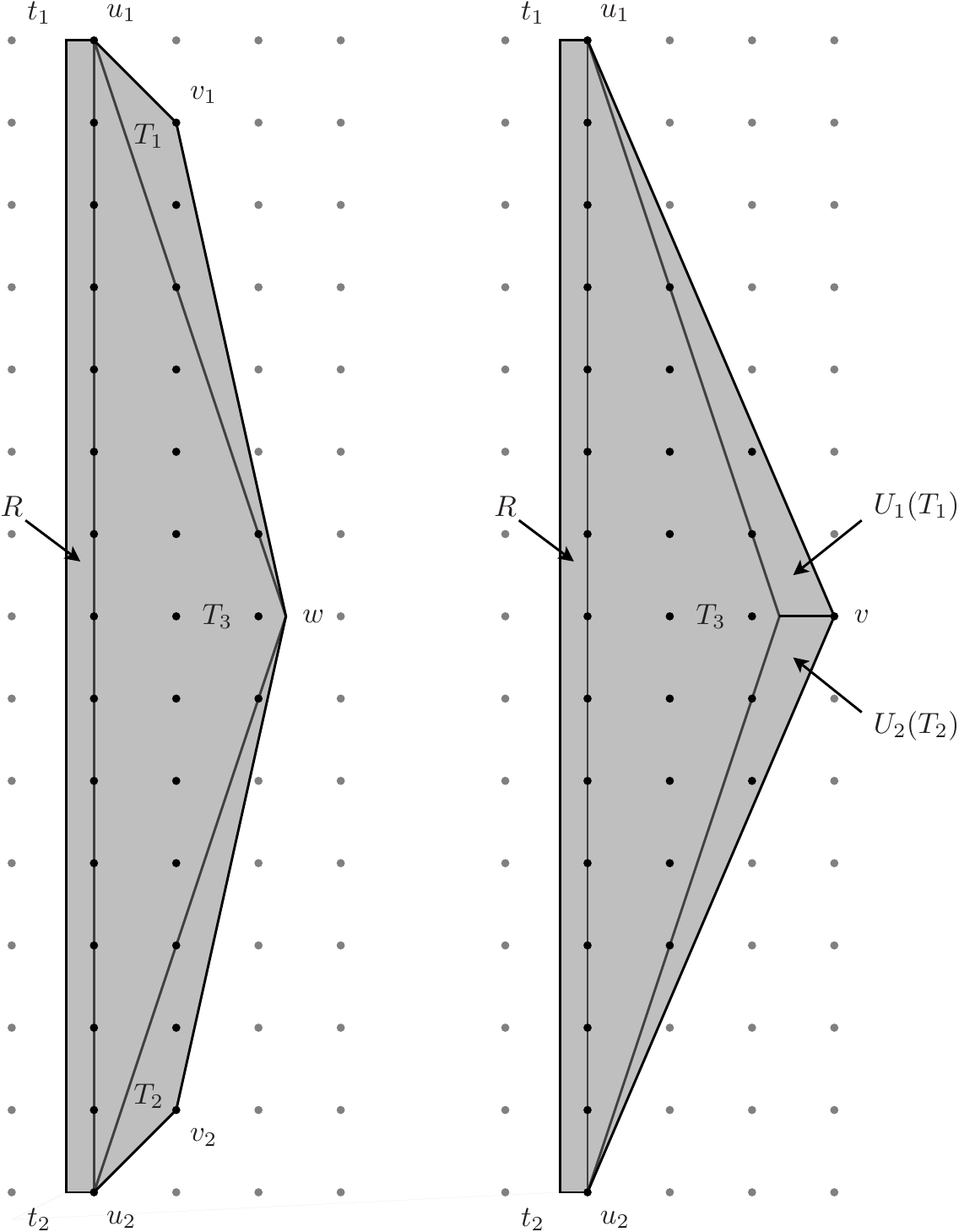}
         \caption{%
            On left: polygon $H$ in the case $s = 3$.  On right:
            polygon $H'$ resulting from unimodular transformation
            of pieces of $H$.
         }%
         \label{fig:CoefficientPeriodsConstruction}
      \end{figure*}
      To show that $H$ has period sequence $(1, s, 1)$, we
      subdivide $H$ into a rectangle and three triangles as
      follows (see left of Figure
      \ref{fig:CoefficientPeriodsConstruction}):
      \begin{align*}
         R     &= \conv\{t_{1}, t_{2}, u_{2}, u_{1}\}, &
         T_{1} &= \conv\{ u_{1}, v_{1}, w \}, \\
         T_{2} &= \conv\{ u_{2}, v_{2}, w \}, &
         T_{3} &= \conv \{ u_{1}, u_{2}, w \}.
      \end{align*}
      Let $v = (s, 0)^{\tran}$.  Write $U_{1} = U_{u_{1} w}^{+}$
      and $U_{2} = U_{u_{2}w}^{-}$.  Then $U_{1}(T_{1}) = \conv
      \setof{u_{1}, v, w}$ and $U_{2}(T_{2}) = \conv \setof{u_{2},
      v, w}$.
      
      Let $H' = R \cup U_{1}(T_{1}) \cup U_{2}(T_{2}) \cup T_{3}$
      (see right of Figure
      \ref{fig:CoefficientPeriodsConstruction}).
      Though $H'$ was formed from unimodular images of pieces of
      $H$, we do not quite have $\Ehr_{H} = \Ehr_{H'}$.  This is
      because each point in the half-open segment $(w, v]$ has two
      pre-images in $H$.  Since this segment is equivalent under a
      unimodular transformation to $h = (\tfrac{1}{s}, 1]$, the
      correct equation is
      \begin{equation}
      \label{eq:EhrHprimeFromEhrH}
	 \Ehr_{H} = \Ehr_{H'} + \Ehr_{h}.
      \end{equation}
              
      Let $T = U_{1}(T_{1}) \cup U_{2}(T_{2}) \cup T_{3}$.  Then
      $T$ is an integral triangle intersecting $R$ along a lattice
      segment, and $H' = R \cup T$.  Hence, by Lemma
      \ref{lem:CoeffsofRectangleUnionIntegral}, $c_{H',1}$ has
      minimum period $s$, and so, by equation
      \eqref{eq:EhrHprimeFromEhrH}, $c_{H,1}$ also has minimum
      period $s$.
      
      It remains only to show that $c_{H,0}$ has minimum period
      $1$.  Again, from \eqref{eq:EhrHprimeFromEhrH}, we have that
      \begin{equation}
      \label{eq:ConstantTerm}
         c_{H,0} = c_{H',0} + c_{h,0}.
      \end{equation}
      From Lemma \ref{lem:CoeffsofRectangleUnionIntegral}, we know
      that $c_{H',0} = c_{\ell,0}$.  Therefore, by
      \eqref{eq:ConstantCoeffofLineSegment}, $c_{H,0}$ is
      identically $1$.
      
      We now construct a triangle with period sequence $(1, 1, t)$
      for integral $t \ge 2$.  Let
      \begin{equation*}
         Q = u_{1} + \conv\setof{(0, 0), (1, -1), (1/t, 0)}.
      \end{equation*}
      McMullen's bound (Theorem \ref{thm:McMullenBound}) implies
      that the minimum period of $c_{Q,1}$ is $1$.  Hence, it
      suffices to show that the minimum quasi-period of $\Ehr_{Q}$
      is $t$.  Observe that $Q$ is equivalent to
      $\conv\setof{(0,0), (1,0), (0, 1/t)}$ under a unimodular
      transformation.  Hence, one easily computes that
      $\sum_{k=0}^{\infty} \Ehr_{Q}(k) \zeta^{k} = (1-\zeta)^{-2}
      (1 - \zeta^{t})^{-1}$.  Note that among the poles of this
      rational generating function are primitive
      $t$\textsuperscript{th} roots of unity.  It follows from the
      standard theory of rational generating functions that
      $\Ehr_{Q}$ has minimum quasi-period $t$ (see, \emph{e.g.},
      \cite[Proposition 4.4.1]{Sta97}).
      
      Finally, for $s, t \ge 2$, let $P = H \cup Q$.  Note that
      $H$ and $Q$ have disjoint interiors, $H \cap Q$ is a lattice
      segment of length $1$, and $H \cup Q$ is convex.  It follows
      that $P$ is a convex polygon and $\Ehr_{P} = \Ehr_{H} +
      \Ehr_{Q} - \Ehr_{[0,1]}$.  Therefore, $P$ has period
      sequence $(1, s, t)$, as required.
   \end{IEEEproof}
   
   
   \def\cprime{$'$}
   
\end{document}